\documentclass[11pt]{article}
\usepackage{amssymb,amsmath,latexsym, verbatim, enumerate}
\usepackage{authblk}
\usepackage[hyperfootnotes=false]{hyperref}
\usepackage[mathscr]{eucal}
\usepackage{mathrsfs}

\usepackage{tikz}
\usetikzlibrary{decorations.fractals}
\usepackage{pgfplots}
\usetikzlibrary{lindenmayersystems}

\allowdisplaybreaks

\begin{document}

\begin{doublespace}

\newtheorem{thm}{Theorem}[section]
\newtheorem{lemma}[thm]{Lemma}
\newtheorem{cond}[thm]{Condition} 
\newtheorem{defn}[thm]{Definition}
\newtheorem{prop}[thm]{Proposition}
\newtheorem{corollary}[thm]{Corollary}
\newtheorem{remark}[thm]{Remark}
\newtheorem{example}[thm]{Example}
\newtheorem{conj}[thm]{Conjecture}
\numberwithin{equation}{section}
\def\ee{\varepsilon}
\def\qed{{\hfill $\Box$ \bigskip}}

\def\b{\mathfrak{b}}
\def\eps{\varepsilon}

\def\j{{\boldsymbol j}}

\def\NN{{\cal N}}
\def\AA{{\cal A}}
\def\MM{{\cal M}}
\def\BB{{\cal B}}
\def\CC{{\cal C}}
\def\LL{{\cal L}}
\def\DD{{\cal D}}
\def\FF{{\cal F}}
\def\EE{{\cal E}}
\def\QQ{{\cal Q}}
\def\RR{{\mathbb R}}
\def\R{{\mathbb R}}
\def\L{{\bf L}}
\def\K{{\bf K}}
\def\S{{\bf S}}
\def\A{{\bf A}}
\def\E{{\mathbb E}}
\def\F{{\bf F}}
\def\P{{\mathbb P}}
\def\N{{\mathbb N}}
\def\eps{\varepsilon}
\def\wh{\widehat}
\def\wt{\widetilde}
\def\pf{\noindent{\bf Proof.} }
\def\beq{\begin{equation}}
\def\eeq{\end{equation}}
\def\lam{\lambda}
\def\H{\mathcal{H}}
\def\nn{\nonumber}
\def\L{\mathcal{L}}
\def\eqd{\,{\buildrel d \over =}\,}
\newcommand{\ud}{\mathrm{d}}

\newcommand\blfootnote[1]{%
  \begingroup
  \renewcommand\thefootnote{}\footnote{#1}%
  \addtocounter{footnote}{-1}%
  \endgroup
}

\title{\Large \bf  Spectral heat content on a class of fractal sets for subordinate killed Brownian motions}
\author{Hyunchul Park and Yimin Xiao}

\date{ \today}
\maketitle


\begin{abstract}
We study the spectral heat content for a class of open sets with fractal boundaries determined by similitudes 
in $\R^{d}$, $d\geq 1$, with respect to subordinate killed Brownian motions via $\alpha/2$-stable subordinators and
establish the asymptotic behavior of the spectral heat content as $t\to 0$ for the full range of $\alpha\in (0,2)$.
Our main theorems show that these asymptotic behaviors depend on whether the sequence of logarithms 
of the coefficients of the similitudes is arithmetic when $\alpha\in [d-\b,2)$, where $\b$ is the 
interior Minkowski dimension of the boundary of the open set. The main tools for proving the theorems 
are the previous results on the spectral heat content for Brownian motions and the renewal 
theorem.
\end{abstract}

\section{Introduction}\label{introduction}

Spectral heat content on an open set
$D\subset \R^{d}$ measures the total heat that remains on $D$ 
at time $t>0$ when the initial temperature is one with Dirichlet boundary condition outside $D$. The 
spectral heat content with respect to Brownian motions has been studied intensively, not only for 
domains with smooth boundary (for example, see \cite{BD1989}) but also for certain domains with fractal 
boundaries such as the $s$-adic von Koch snowflake (see \cite{B1994, B2000, BG1998, FLV1995}).

Recently there have been increasing interests in the spectral heat content for more general L\'evy 
processes (see \cite{V1, GPS19, KP, PS}). In \cite{GPS19} the authors studied the spectral heat 
content for some class of L\'evy processes on bounded open sets in $\R$. In particular, it is proved in 
\cite[Theorems 4.2 and 4.14]{GPS19} that when the underlying open set has infinitely many components 
(or infinitely many non-adjacent components in the case of the Cauchy process) the decay rate of the
spectral heat content  is strictly bigger than that of the spectral heat content with 
respect to open sets with finitely many components. Hence, a natural question is to determine the exact decay rate 
of the spectral heat content for L\'evy processes when there are infinitely many components in $D\subset \R$. 

Many L\'evy processes can be realized as subordinate (time-changed) Brownian motions, 
where the time change is given by an independent subordinator. Since we need two operations to define 
the spectral heat content for subordinate Brownian motions, time-change and killing, there are two objects 
that can be called the spectral heat content for subordinate Brownian motions; One is related to the 
{\it killed subordinate Brownian motions} (do time-change first then kill the time-changed Brownian motions 
when they exit the domain under consideration) and the other is related to {\it subordinate killed Brownian motions}
 (kill Brownian motions when they exit the domain, then make time-change for the killed Brownian motions). 
Even though the spectral heat content for killed subordinate Brownian motions is a natural object to study as 
it covers a large class of the spectral heat content for killed L\'evy processes, the spectral heat content for 
subordinate killed Brownian motions is also important as it oftentimes gives useful information on the spectral 
heat content for the killed subordinate Brownian motions. For example, in \cite{PS21} the asymptotic behavior of the 
spectral heat content for subordinate killed Brownian motions with respect to stable subordinators provides 
crucial information on the spectral heat content for killed stable processes.

In this paper, we study the spectral heat content for  the subordinate killed 
Brownian motions  (see (\ref{Def:hc}) below for the definition)
when the underlying subordinator is a stable subordinator $S^{(\alpha/2)}=\{S_{t}^{(\alpha/2)}\}_{t\geq 0}$
whose Laplace transform is given by 
\begin{equation}\label{Eq:Laplace}
\E[e^{-\lam S_{t}^{(\alpha/2)}}]=e^{-t\lam^{\alpha/2}}, \quad \lam>0,\ \alpha \in (0,2),
\end{equation}
and the underlying open sets have fractal boundaries which are determined by similitudes in $\R^{d}$. 
Our main results answer the aforementioned question and they show that  the exact decay rates of the 
spectral heat content depend on whether the sequence of logarithms of the coefficients of the similitudes 
is arithmetic when $\alpha\in [d-\b,2)$, where $\b\in (d-1,d)$ is the interior Minkowski 
dimension of the boundary of the underlying set (cf. (\ref{Def:Mdim}) below for definition). It is noteworthy 
to observe that when $\alpha\in (0,d-\b)$ the theorem is independent of whether the sequence $\{\ln(1/r_{j})\}$
 of the logarithms of the coefficients of the similitudes is arithmetic or not.

The main technique for studying the asymptotic behavior of the spectral heat content for open sets
with fractal boundaries in this paper is \textit{the renewal theorem} in \cite{LV1996}. Two crucial properties 
that we need are the additivity property of the spectral heat content under disjoint union (Lemma 
\ref{lemma:disjoint sum}) and the scaling property of the subordinate killed Brownian motions with 
respect to stable subordinators (Lemma \ref{lemma:scaling}). However, in order to apply the renewal 
theorem, one needs an exponential decay condition \eqref{cond:RT2} and this is only valid 
when $\alpha\in(d-\b,2)$. In order to establish the asymptotic behavior of the spectral heat content as $t\to 0$ 
for the full range of $\alpha\in (0,2)$, 
we will use the weak convergence of L\'evy measure (see Proposition \ref{prop:alpha small}) to establish 
the asymptotic behavior for the case when $\alpha\in (0,d-\b)$. The result is given as Theorem \ref{thm:small alpha}. 
The remaining case 
when $\alpha=d-\b\in(0,1)$  is proved in Theorems \ref{thm:d-b arithmetic} and \ref{thm:d-b non-arithmetic}
for the arithmetic and non-arithmetic cases, respectively. 
We observe that there is an extra logarithm term $\ln(1/t)$ in the decay rate of the spectral heat content 
when $\alpha=d-\b$. This is due to the fact that the heat loss $|G|-Q_{G}^{(2)}(u)$ for 
Brownian motions on the open set $G$ with fractal boundaries, where $|G|-Q_{G}^{(2)}(u)=\int_{G}\P_{x}(\tau_{G}^{(2)} 
\leq u)dx$ and $\tau_{G}^{(2)}$ is the first exit time of the Brownian motions out of $G$, is barely 
non-integrable with respect to the law of stable subordinator $S_{t}^{(\alpha/2)}$.
Note that this occurrence of extra logarithm term $\ln(1/t)$ is observed for smooth open sets in \cite{Val2016, PS,PS21}, 
but it happens at a different index $\alpha$. More specifically, for the spectral heat content on smooth sets this phenomenon 
happens when $\alpha=1$,  whereas for open sets with fractal boundaries as in our case, this happens when $\alpha=d-\b$, 
which is strictly less than 1.

The organization of this paper is as follows. In Section \ref{section:preliminaries} we set up notations, define the class of
open sets with fractal boundaries in \eqref{eqn:SS}, and recall some facts which will be used for proving our main theorems. 
In particular, we recall the renewal theorem in \cite{LV1996} and the result for the spectral heat content for Brownian motions 
in Theorem \ref{thm:BM}. Section \ref{section:SHC} is the main part of this paper and here we study the spectral heat 
content for subordinate killed Brownian motions. The main results are Theorems \ref{thm:large alpha}, \ref{thm:d-b arithmetic},
 \ref{thm:d-b non-arithmetic}, and \ref{thm:small alpha}. This section is divided into three 
subsections for the cases $\alpha\in (d-\b,2)$, $\alpha=d-\b$, and $\alpha\in (0,d-\b)$, respectively.

We use $c_i$ to denote constants whose values are unimportant and may change from one appearance to another.
The notations $\P_{x}$ and $\E_{x}$ mean probability and expectation of the underlying processes started at $x\in \R^{d}$, 
and we use $\P=\P_{0}$ and $\E=\E_{0}$ to simplify notations.

\section{Preliminaries}\label{section:preliminaries}

In this section, we introduce notations and recall some facts  that will be used for proving the main theorems in Section 3. 

\subsection{Some geometric notions}
We first recall some geometric notions and the definition of the class of open sets with fractal boundaries from \cite{LV1996}. 
See also \cite{Lapidus13, Lapidus17} for more recent developments.

For any bounded (open) set $G \subset \R^d$ with boundary $\partial G$ and any $\varepsilon > 0$, let 
\[
G_\eps^{\rm int} = \{x \in G: {\rm dist}(x, \partial G) < \eps\}
\]
be the interior Minkowski sausage of radius $\eps$ of the boundary $\partial G$. We denote by $\mu (\eps; G)$ 
the $d$-dimensional Lebesgue measure of $G_\eps^{\rm int} $. For any $ s >0$, define
\[
{\mathscr M }^*(s, \partial G) = \limsup_{\eps \to 0} {\eps}^{-(d-s)} \mu (\eps; G)
\]
and
\[
{\mathscr M}_* (s, \partial G) = \liminf_{\eps \to 0} {\eps}^{-(d-s)} \mu (\eps; G).
\]
Following  \cite[Definition 1.2]{LV1996}, the interior Minkowski dimension of $\partial G$ 
(which is also called the Minkowski dimension of $\partial G$ relative to $G$)  is defined by
\begin{equation} \label{Def:Mdim}
\dim_{\rm M}^{\rm int} (\partial G)  = \inf\{s > 0: \, {\mathscr M }^*(s, \partial G) = 0 \} 
= \sup\{ s > 0: \, {\mathscr M }_*(s, \partial G) = \infty \}.
\end{equation}
If, for $s = \dim_{\rm M}^{\rm int} (\partial G)$,  we have $0< {\mathscr M }^*(s, \partial G) 
= {\mathscr M}_* (s, \partial G) < \infty$, then, $\partial G$  is said to be Minkowski measurable 
relative to $G$.

\begin{defn}
A map $R:\R^{d}\to \R^{d}$ is called a \textit{similitude} with coefficient $r>0$ if 
$$
|Rx-Ry|=r|x-y| \quad \hbox{ for all } x,\, y\in \R^{d}.
$$
\end{defn}
It is well known (cf. e.g., \cite[p.191]{LV1996}) that any similitude is a composition of a homothety with 
coefficient $r$, an orthogonal transform, and a translation. 

Now we define the class of open sets with fractal boundaries that we will consider in this paper. 
Let $G_{0} \subset \R^{d}$ be a bounded open set. When $d=1$ we assume that $G_0$ is a bounded
open interval and when $d\geq 2$ we assume that $G_{0}$ is a bounded $C^{1,1}$ open set. Let 
$R_{j}$ ($1\leq j\leq N$) be similitudes with coefficients $r_{j}$, respectively. 

For each $n \ge 1$, define $\Upsilon_n = \{ {\boldsymbol j} = (j_1, \ldots, j_n),\ 1 \le j_i \le N\}$.  
We define the set $G$ by 
\beq\label{eqn:SS1}
G= \bigg(\bigcup_{n=1}^{\infty} \bigcup_{ {\boldsymbol j } \in \Upsilon_n}  {\mathscr R }_{\j}G_0 \bigg)\cup G_{0},
\eeq
where,  for every $\j = (j_1, \ldots, j_n) \in \Upsilon_n$, ${\mathscr R }_{\j}$ is the similitude  defined by ${\mathscr R }_{\j} 
=R_{j_1}\circ \cdots \circ R_{j_n}$. 
It follows from (\ref{eqn:SS1}) that $G$ can be represented as
\beq\label{eqn:SS}
G= \bigg(\bigcup_{j=1}^{N}R_{j}G\bigg) \cup G_{0}.
\eeq
We assume that all the sets $R_{j}G $, $1 \le j \le N$, and $G_0$ in \eqref{eqn:SS} are pairwise disjoint. 
As in \cite[Equation (1.7)]{LV1996} we also assume that $\sum_{j=1}^{N}r_{j}^{d}<1<\sum_{j=1}^{N}r_{j}^{d-1}$.
Since the expressions in \eqref{eqn:SS} are pairwise disjoint, we have
\beq\label{eqn:measure}
|G|=\sum_{j=1}^{N}r_{j}^{d}|G|+|G_{0}|.
\eeq
Also, the condition $\sum_{j=1}^{N}r_{j}^{d}<1<\sum_{j=1}^{N}r_{j}^{d-1}$ ensures that $G$ has a finite volume 
and there exists a unique number $\b\in (d-1,d)$ such that 
\beq\label{eqn:MD}
\sum_{j=1}^{N}r_{j}^{\b}=1.
\eeq
It follows from \cite[Theorem A]{LV1996} that the number $\b$ is equal to the interior Minkowski dimension of $\partial G$. 

As an illustration, we consider the following examples. Let $R_1$ and $R_2$ be two similitudes on $\R$ defined by
$$
R_1(x) = \frac 1 3 x \quad \hbox{ and } \quad R_2(x) = \frac 1 3 x + \frac 2 3.
$$
We take $G_0 = (\frac13, \frac23)$. Then, it is easy to observe that the set defined in \eqref{eqn:SS1} is given by
 $G=(0,1)\setminus \mathfrak{C}$, where $\mathfrak{C}$ is the standard ternary Cantor set and $G$ satisfies \eqref{eqn:SS},
and $\b = \log 2/\log 3$. Similarly, the open set $G \subset (0, \infty)^2$ with the Sierp\'inski gasket as its 
boundary can be obtained from  $G_0 $ being the open triangle with vertices $(1/4, \sqrt{3}/ 4),\, (1/2, 0)$ and $(3/4, \, \sqrt{3}/ 4)$ 
(notice that the boundary $\partial G_0$ is not $C^{1,1}$) and  three similitudes on $\R^2$ defined by
\[
R_1(x) = \frac 1 2 x, \quad R_2(x) = \frac 1 2 x + \Big( \frac 1 4, \frac {\sqrt{3} } 4 \Big)  \ \ \  \hbox{ and } \ \  R_3(x) = \frac 1 2 x +
\Big ( \frac 1 2, 0\Big).
\]
The (interior) Minkowski dimension of the boundary $\partial G$ is $\b = \log 3/\log 2$.

\subsection{The renewal theorem}

Now we state a version of the renewal theorem from \cite{LV1996}.
Let $f:\R\to \R$ be a map. 
For any $\gamma\in \R$ define 
$$
L_{\gamma}f(z)=f(z-\gamma),
$$
and 
$$
Lf(z)=\sum_{j=1}^{N}c_{j}L_{\gamma_{j}}f(z)=\sum_{j=1}^{N}c_{j}f(z-\gamma_{j}),
$$
where $c_{j}>0$, $\gamma_{j}$ are distinct points in $\R$, and $\sum_{j=1}^{N}c_{j}=1$.

Consider the following \textit{renewal equation}
\beq\label{eqn:RE}
f=Lf +\phi.
\eeq
Intuitively, it is natural to expect that the solution of the renewal equation is given by
\beq\label{eqn:sol}
f(z)=\sum_{n=0}^{\infty}L^{n}\phi(z)=\phi(z)+\sum_{n=1}^{\infty}\sum_{c_{i_1},\cdots, c_{i_n}} 
c_{i_1}\cdots c_{i_n}L_{\gamma_{i_1}}\cdots L_{\gamma_{i_{n}}}\phi(z).
\eeq
The following renewal theorem says it is indeed the case under certain conditions.
We say a set of finite real numbers $\{\gamma_{1},\cdots, \gamma_{N}\}$ is \textit{arithmetic}
 if $\frac{\gamma_{i}}{\gamma_{j}}\in \mathbb{Q}$ for all indices. The maximal number $\gamma$ 
 such that $\frac{\gamma_{i}}{\gamma}\in \mathbb{Z}$ is called the span of $\{\gamma_{1},\cdots, \gamma_{N}\}$.
If the set is not arithmetic, it is called \textit{non-arithmetic}.
\begin{thm}[Renewal Theorem \cite{LV1996}]\label{thm:RT}
Suppose that a map $f:\R\to \R$ satisfies the renewal equation \eqref{eqn:RE} and it satisfies
\beq\label{cond:RT1}
\lim_{z\to-\infty}f(z)=0,
\eeq
and 
\beq\label{cond:RT2}
|\phi(z)|\leq c_{1}e^{-c_{2}|z|}, \quad z\in \R,
\eeq
for some constants $c_1,c_2>0$. Then, the solution of the renewal equation \eqref{eqn:RE} is given by \eqref{eqn:sol}. 
Furthermore, if $\{\gamma_{j}\}$ is non-arithmetic, then 
$$
f(z)=\frac{1}{\sum_{j=1}^{N}c_{j}\gamma_{j}}\int_{-\infty}^{\infty}\phi(x)dx +o(1), \text{ as } z\to\infty.
$$ 
If $\{\gamma_{j}\}$ is arithmetic with span $\gamma$, then
$$
f(z)=\frac{\gamma}{\sum_{j=1}^{N}c_{j}\gamma_{j}}\sum_{k=-\infty}^{\infty}\phi(z-k\gamma) +o(1), \text{ as } z\to\infty.
$$
\end{thm}

\subsection{The spectral heat content of subordinate killed Brownian motions}

For an open set $G\subset \R^{d}$ we define the spectral heat content for Brownian motion
$W=\{W_{t}\}_{t\geq	 0}$ on $G$ as
$$
Q_{G}^{(2)}(t)=\int_{G}\P_{x} \big(\tau_{G}^{(2)} > t \big)dx, \quad \tau_{G}^{(2)}=\inf\{t>0 : W_{t}\notin G\}.
$$

We record the following lemma from \cite[Lemma 4.4]{LV1996}. 
Note that there was a typo there and $r^{2}$ should be written as $r^{d}$. 
\begin{lemma}[Lemma 4.4 \cite{LV1996}]\label{lemma:scaling BM}
Let $G$ be an open set in $\R^{d}$ and $R$ be a similitude with coefficient $r$. Then, 
$$
Q_{RG}^{(2)}(t)=r^{d}Q_{G}^{(2)}(t/r^2).
$$
\end{lemma}

We recall the following theorem for the spectral heat content for Brownian motion 
from \cite[Theorem D]{LV1996}.
\begin{thm}[Theorem D \cite{LV1996}]\label{thm:BM}
Let $G$ be a set defined as in \eqref{eqn:SS1}
with $R_{j}$ being similitude with coefficient $r_{j}$, 
and $G_{0}$ is either a bounded open interval when $d=1$, or a bounded $C^{1,1}$ open set when $d\geq 2$. 
\begin{enumerate}
\item 
If $\{\ln(\frac{1}{r_{j}})\}_{j=1}^{N} $ is non-arithmetic, then  
$$
Q_{G}^{(2)}(t)=|G|-Ct^{\frac{d-\b}{2}} + o(t^{\frac{d-\b}{2}}),\quad \text{ where }
C=\frac{\int_{0}^{\infty}\left(|G_{0}|-Q_{G_{0}}^{(2)}(u)\right)u^{-(1+\frac{d-\b}{2})}du}{\sum_{j=1}^{N}r_{j}^{\b}\ln(\frac{1}{r_{j}^{2}})}.
$$
\item 
If $\{\ln(\frac{1}{r_{j}})\}_{j=1}^{N}$ is arithmetic with span $\rho$,  then
$$
Q_{G}^{(2)}(t)=|G|-s(-\ln t)t^{\frac{d-\b}{2}}+ o(t^{\frac{d-\b}{2}}), \quad \text{ as } t\to 0,
$$
where the function $s(\cdot)$ is defined by 
\begin{eqnarray*}
s(z)&:=&\frac{2\rho}{\sum_{j=1}^{N}(r_{j})^{\b}\ln(\frac{1}{r_{j}^{2}})}\sum_{n=-\infty}^{\infty}\left(|G_{0}|-
Q_{G_{0}}^{(2)}\big(e^{-(z-2n\rho)}\big)\right)e^{\frac{d-\b}{2}(z-2n\rho)}.
\end{eqnarray*}
\end{enumerate}
\end{thm}

Now we introduce the spectral heat content for subordinate killed Brownian motions. 
Let $W=\{W_{t}\}_{t\geq 0}$ be Brownian motion in $\R^{d}$ and let 
$S^{(\alpha/2)} =\{ S_{t}^{(\alpha/2)} \}_{t\geq 0}$ be an $(\alpha/{2})$-stable 
subordinator with Laplace transform given by (\ref{Eq:Laplace}),
which is independent of $W$. 
Let $D$ be any open set in $\R^{d}$. Then, the spectral heat content $\tilde{Q}_{D}^{(\alpha)}(t)$ for 
subordinate killed Brownian motions with respect to stable subordinator $S^{(\alpha/2)}$ on $D$ is defined by
\begin{equation}\label{Def:hc}
\tilde{Q}_{D}^{(\alpha)}(t)=\int_{D}\P_{x}\Big(\tau_{D}^{(2)} > S_{t}^{(\alpha/2)}\Big)dx,
\end{equation}
where $\tau_{D}^{(2)} = \inf\{t>0 : W_{t}\notin D\}$.

We will need the following important properties for $\tilde{Q}_{D}^{(\alpha)}(t)$; One is the additivity 
under disjoint union and another is the scaling property.
\begin{lemma}\label{lemma:disjoint sum}
Let $D_{1},D_{2}$ be open sets in $\R^{d}$ with $D_{1}\cap D_{2}=\emptyset$. Then, 
$$
\tilde{Q}^{(\alpha)}_{D_1\cup D_{2}}(t)=\tilde{Q}^{(\alpha)}_{D_1}(t)+\tilde{Q}^{(\alpha)}_{D_2}(t).
$$
\end{lemma}
\pf
Note that
\begin{equation*}
\begin{split}
\tilde{Q}^{(\alpha)}_{D_1\cup D_{2}}(t) &=\int_{D_{1}\cup D_{2}}\P_{x}\Big(\tau^{(2)}_{D_{1}\cup D_{2}} >S_{t}^{(\alpha/2)}\Big)dx\\
&=\int_{D_{1}}\P_{x}\Big(\tau^{(2)}_{D_{1}\cup D_{2}} >S_{t}^{(\alpha/2)}\Big)dx 
+\int_{D_{2}}\P_{x}\Big(\tau^{(2)}_{D_{1}\cup D_{2}} >S_{t}^{(\alpha/2)}\Big)dx.
\end{split}
\end{equation*}
Note that under $\P_{x}$ with $x\in D_{1}$ we have
$$
\tau^{(2)}_{D_{1}\cup D_{2}}=\inf\{t>0 : W_{t}\notin D_{1}\cup D_{2}\}=\inf\{t>0 : W_{t}\notin D_{1}\}=\tau_{D_{1}}^{(2)}.
$$
Hence, we have
$$
\int_{D_{1}}\P_{x}\Big(\tau^{(2)}_{D_{1}\cup D_{2}} >S_{t}^{(\alpha/2)}\Big)dx
=\int_{D_{1}}\P_{x}\Big(\tau^{(2)}_{D_{1}} >S_{t}^{(\alpha/2)}\Big)dx=\tilde{Q}_{D_{1}}^{(\alpha)}(t).
$$
It can be proved in the same way that the integral on $D_2$ gives $\tilde{Q}^{(\alpha)}_{D_2}(t).$
\qed

\begin{remark}
Let $Q_{D}^{(\alpha)}(t):=\int_{D}\P_{x}\Big(\tau_{D}^{(\alpha)}>t \Big)dx$ be the spectral heat content for killed stable processes, 
where $\tau_{D}^{(\alpha)}$ is the first exit time of the $\alpha$-stable process $W_{S^{(\alpha/2)}}=\{W_{S_{t}^{(\alpha/2)}}\}_{t\geq 0}$. This is the 
spectral heat content related to the killed subordinate Brownian motions by stable subordinators. 
Note that for disjoint sets $D_{1}$ and $D_{2}$ we have 
\begin{align*}
Q^{(\alpha)}_{D_1 \cup D_{2}}(t)&=
\int_{D_1\cup D_{2}}\P\Big(\tau^{(\alpha)}_{D_1 \cup D_{2}}>t\Big)dx 
= \int_{D_1}\P\Big(\tau^{(\alpha)}_{D_1 \cup D_{2}}>t \Big)dx 
+\int_{ D_{2}}\P\Big(\tau^{(\alpha)}_{D_1 \cup D_{2}}>t \Big)dx \\
&\geq  \int_{D_1}\P\Big(\tau^{(\alpha)}_{D_1 }>t\Big)dx 
+\int_{ D_{2}}\P\Big(\tau^{(\alpha)}_{ D_{2}}>t \Big)dx 
=Q^{(\alpha)}_{D_1}(t)+Q^{(\alpha)}_{D_{2}}(t).
\end{align*}
Furthermore, the inequality can be strict as $\tau^{(\alpha)}_{D_1 \cup D_{2}} \neq \tau^{(\alpha)}_{D_1}$ when the process starts 
at $x\in D_{1}$ because the process starting at $x\in D_{1}$ can jump into $D_{2}$ without visiting the complement of $D_{1}\cup D_{2}$.
Hence, the spectral heat content for killed subordinate Brownian motions does not satisfy the additivity property under disjoint union.
\end{remark}

\begin{lemma}\label{lemma:scaling}
Let $R$ be a similitude with coefficient $r$ and $G$ is any open set in $\R^{d}$. Then, we have
$$
\tilde{Q}^{(\alpha)}_{RG}(t)=r^{d}\tilde{Q}^{(\alpha)}_{G}(t/r^{\alpha}), \quad t>0.
$$
\end{lemma}
\pf
By the scaling property and rotational invariance of Brownian motions, we observe that 
$$
\tau_{RG}^{(2)} \text{ under } \P_{Rx} \text{ is equal in distribution to } r^{2}\tau_{G}^{(2)} \text{ under } \P_{x}.
$$
By the change of variable $x=Ry$ and the scaling property of $S_{t}^{(\alpha/2)}$ we have
\begin{equation*}
\begin{split}
\tilde{Q}^{(\alpha)}_{RG}(t)&=\int_{RG}\P_{x} \Big(\tau_{RG}^{(2)}>S_{t}^{(\alpha/2)}\Big)dx
=\int_{G}\P_{Ry} \Big(\tau_{RG}^{(2)}>S_{t}^{(\alpha/2)} \Big)r^{d}dy\\
&=\int_{G}\P_{y}\Big(r^{2}\tau_{G}^{(2)}>S_{t}^{(\alpha/2)} \Big)r^{d}dy
=\int_{G}\P_{y} \Big(\tau_{G}^{(2)}>r^{-2}S_{t}^{(\alpha/2)}\Big)r^{d}dy\\
&=\int_{G}\P_{y}\Big(\tau_{G}^{(2)}>S_{tr^{-\alpha}}^{(\alpha/2)}\Big)r^{d}dy
=r^{d}\tilde{Q}_{G}^{(\alpha)}(t/r^{\alpha}).
\end{split}
\end{equation*}
This proves the lemma. 
\qed

\section{Asymptotic behavior of the spectral heat content}\label{section:SHC}

\subsection{The case of $\alpha \in (d- \b,2)$}
Analogous to \cite[Theorem D]{LV1996}, we will prove that the spectral heat content 
$\tilde{Q}^{(\alpha)}_{G}(t)$ as defined in (\ref{Def:hc}) has the form 
$$
\tilde{Q}^{(\alpha)}_{G}(t)=|G|-f(-\ln t)t^{\frac{d- \b}{\alpha}} +o(t^{\frac{d- \b}{\alpha}}),
$$
when $\alpha\in (d-\b,2)$, where $\b$ is the constant in \eqref{eqn:MD}.
We start with the following lemma. 

\begin{lemma}\label{lemma:condition for RT}
Assume that $\alpha \in (d-\b,2)$.
Suppose that $G_{0}$ is an open interval when $d=1$ or a bounded $C^{1,1}$ open set when $d\geq 2$. 
Define
$$
\psi(z)=\left(|G_{0}|-\tilde{Q}_{G_{0}}^{(\alpha)}(e^{-z})\right)e^{\frac{d-\b}{\alpha}z}.
$$
Then, there exists a constant $c=c(\alpha,\b,d)>0$ such that 
$$
|\psi(z)|\leq ce^{-c|z|} \text{ for all }z\in \R.
$$
\end{lemma}
\pf
We define  the function
$$
\phi(t)=\left(|G_{0}|-\tilde{Q}^{(\alpha)}_{G_{0}}(t) \right)t^{-\frac{d-\b}{\alpha}}, \ \ \ t > 0,  
$$
so that $\psi(z)=\phi(e^{-z})$.

The case when $z\to -\infty$ is easy since we have 
$$
|\phi(t)|\leq |G_{0}|t^{-\frac{d-\b}{\alpha}}
$$
and this gives
$$
\psi(z)=\phi(e^{-z})\leq |G_{0}|e^{\frac{d- \b}{\alpha}z}=|G_{0}|e^{-\frac{d-\b}{\alpha}|z|}
$$
for all $z\leq 0$ and $\alpha \in (0,2)$.

Now we handle the case when $z\to \infty$, or $t = e^{-z}\to 0$.
First, assume that $\alpha \in (1,2)$.
Since $G_{0}$ is an interval when $d=1$ or a bounded $C^{1,1}$ open set when $d\geq 2$, it follows 
from \cite[Theorem 1.1]{PS} that 
there exists a constant $c_1>0$ such that
$$
|G_{0}|-\tilde{Q}^{(\alpha)}_{G_{0}}(t)\leq c_{1}t^{1/\alpha}
$$
for all $0<t\leq 1$.
Hence, 
$$
\phi(t)\leq c_{1}t^{\frac{( \b+1)-d}{\alpha}} \text{ for } 0<t\leq 1.
$$
Since $\b\in (d-1,d)$ we note that $\frac{(\b+1)-d}{\alpha}>0$.
Let $z=-\ln t$ and we conclude that 
$$
\phi(t)=\phi(e^{-z})=\psi(z)\leq c_{2}e^{-c_{3}|z|}, \quad z\in \R,
$$
where $c_{2}=\max(c_1, |D_0|)$ and $c_{3}=\min(\frac{\b+1-d}{\alpha}, \frac{d-\b}{\alpha})>0$.

Second, when $\alpha=1$ we have from \cite[Theorem 1.1]{PS}
$$
|G_{0}|-\tilde{Q}^{(\alpha)}_{G_{0}}(t)\leq c_{2}t\ln(1/t)
$$
for all $0<t\leq 1$.
Hence, 
$$
\phi(t)\leq c_{2}t^{1-(d-\b)}\ln(1/t)=c_{2}t^{(\b+1)-d}\ln(1/t) \quad \text{ for } 0<t\leq 1,
$$
and this implies
$$
\psi(z)=\phi(e^{-z})\leq c_{2}ze^{-z(\b+1-d)} \text{ for } z\geq 0.
$$
Since $\b+1-d>0$ there exists $c_{4}$ and $\eta>0$ such that 
$$
\psi(z)\leq c_{4}e^{-\eta z}   \ \ \text{ for all } z\geq 0.
$$

Finally, we handle the case when $\alpha \in (d-\b,1)$. From \cite[Theorem 1.1]{PS} we have
$$
|G_{0}|-\tilde{Q}^{(\alpha)}_{G_{0}}(t)\leq c_{5}t  
$$
for all $0<t\leq 1$, and this implies $\phi(t)\leq c_{5}t^{\frac{\alpha-(d-\b)}{\alpha}}   \text{ for } t\leq 1$,
which in turn implies $\psi(z)\leq c_{5}e^{-\frac{\alpha-(d-\b)}{\alpha}z} \text{ for } z\geq 0$.
\qed

Here is the main theorem for the case of $\alpha \in (d- \b,2)$.
\begin{thm}\label{thm:large alpha}
Let $\alpha\in (d-\b,2)$, where $\b$ is the constant in \eqref{eqn:MD} and $G$ is a set given as \eqref{eqn:SS1}
with $G_{0}$ being an open interval when 
$d=1$ or a bounded $C^{1,1}$ open set when $d\geq 2$.
If $\{\ln\frac{1}{r_{j}}\}_{j=1}^{N}$ is non-arithmetic, then we have
$$
\tilde{Q}^{(\alpha)}_{G}(t)=|G|-C_{1}t^{\frac{d-\b}{\alpha}}+o(t^{\frac{d-\b}{\alpha}}) \text{ as } t\to 0,
$$
where
\begin{eqnarray*}
C_{1}&=&\frac{\int_{-\infty}^{\infty}\left(|G_{0}|-\tilde{Q}_{G_{0}}^{(\alpha)}(e^{-z})\right)e^{\frac{(d-\b)z}
{\alpha}}dz}{\sum_{j=1}^{N}r_{j}^{\b}\ln(1/r_{j}^{\alpha})}
=\frac{\int_{0}^{\infty}\left(|G_{0}|-\tilde{Q}_{G_0}^{(\alpha)}(t)\right)t^{-1-\frac{d-\b}{\alpha}}dt}
{\sum_{j=1}^{N}r_{j}^{\b}\ln(1/r_{j}^{\alpha})}.
\end{eqnarray*}
If $\{\ln\frac{1}{r_{j}}\}_{j=1}^{N}$ is arithmetic with span $\rho$, then we have
$$
\tilde{Q}^{(\alpha)}_{G}(t)=|G|-f(-\ln t)t^{\frac{d-\b}{\alpha}} +o(t^{\frac{d-\b}{\alpha}}) \text{ as } t\to0,
$$
where
$$
f(z)=\frac{\alpha\rho}{\sum_{j=1}^{N}r_{j}^{\b}\ln(1/r_{j}^{\alpha})}\sum_{n=-\infty}^{\infty}\left(|G_{0}|-\tilde{Q}_{G_{0}}^{(\alpha)}
(e^{-(z-\alpha n\rho)})\right)e^{\frac{d-\b}{\alpha}(z-\alpha n \rho)}.
$$
\end{thm}
\pf
We set 
\beq\label{eqn:expression}
\tilde{Q}_{G}^{(\alpha)}(t)=|G|-f(-\ln t)t^{\frac{d-\b}{\alpha}},
\eeq
and will show that $f(\cdot)$ satisfies the renewal equation and conditions for the renewal theorem. 

Since $G=\bigcup_{j=1}^{N}R_{j}G\cup G_{0}$ and all expressions are disjoint, it follows from Lemmas 
\ref{lemma:disjoint sum} and \ref{lemma:scaling}
\begin{eqnarray*}
&&\tilde{Q}_{G}^{(\alpha)}(t)=\tilde{Q}_{\bigcup_{j=1}^{N}R_{j}G\cup G_{0}}^{(\alpha)}(t)
=\sum_{j=1}^{N}\tilde{Q}_{R_{j}G}^{(\alpha)}(t) +\tilde{Q}_{ G_{0}}^{(\alpha)}(t)=\sum_{j=1}^{N}r_{j}^{d}
\tilde{Q}_{G}^{(\alpha)}(t/r_{j}^{\alpha}) +\tilde{Q}^{(\alpha)}_{G_{0}}(t).
\end{eqnarray*}
Using \eqref{eqn:expression} we have
\[
\begin{split}
&|G|-f(-\ln t)t^{\frac{d-b}{\alpha}}
=\sum_{j=1}^{N}r_{j}^{d}\left(|G|-f(-\ln\frac{t}{(r_{j})^{\alpha}})(\frac{t}{(r_{j})^{\alpha}})^{\frac{d-\b}{\alpha}}\right) +\tilde{Q}^{(\alpha)}_{G_{0}}(t)\\
&=\sum_{j=1}^{N}r_{j}^{d}|G|-\sum_{j=1}^{N}r_{j}^{\b}\cdot f\left(-\ln t -\ln (\frac{1}{(r_{j})^{\alpha}})\right)t^{\frac{d-\b}{\alpha}}
+\left(|G_{0}|-(|G_{0}|-\tilde{Q}^{(\alpha)}_{G_{0}}(t) )\right).
\end{split}
\]
Using \eqref{eqn:measure} we conclude that 
$$
f(-\ln t)=\sum_{j=1}^{N}r_{j}^{\b}\cdot f\left(-\ln t -\ln (\frac{1}{(r_{j})^{\alpha}})\right) +\phi(t), \quad
\text{where } 
\phi(t)=\left(|G_{0}|-\tilde{Q}^{(\alpha)}_{G_{0}}(t) \right)t^{-\frac{d-\b}{\alpha}}.
$$
By changing the variable $z=-\ln t$ we have 
$$
f(z)=\sum_{j=1}^{N}r_{j}^{\b}\cdot f\left(z -\ln (\frac{1}{(r_{j})^{\alpha}})\right) +\phi(e^{-z}).
$$

Note that from \eqref{eqn:expression} we have 
$$
\lim_{z\to-\infty}f(z)=\lim_{t\to \infty}f(-\ln t)=\lim_{t\to \infty}\left(|G|-\tilde{Q}^{(\alpha)}_{G}(t)\right)t^{-\frac{d-\b}{\alpha}}=0,
$$
and this shows that the condition \eqref{cond:RT1} holds.
It follows from Lemma \ref{lemma:condition for RT} that for any $\alpha\in (d-\b,2)$ there exist two constants 
$c_{1},c_{2}>0$ such that
$$
\psi(z)=\phi(e^{-z})\leq c_{1}e^{-c_{2}|z|}\ \  \text{ for all } z\in \R,
$$
and the condition \eqref{cond:RT2} holds. 
Now the conclusions of the theorem follow immediately from the Renewal Theorem \ref{thm:RT}.
\qed

\subsection{The case of $\alpha=d-\b$}

In this subsection, we study the case when $\alpha=d-\b \in (0,1)$.
We need a simple lemma which is similar to \cite[Lemma 3.2]{PS}. The proof is essentially the same with 
obvious modifications and will be omitted.
\begin{lemma}\label{lemma:SS expectation}
For any $\delta>0$ and $\alpha\in(0,2)$, we have
$$
\lim_{t\to 0}\frac{\E \left[\big(S_{1}^{(\alpha/2)}\big)^{\alpha/2}, 0<S_{1}^{(\alpha/2)}<\delta t^{-2/\alpha} \right]}{\ln(1/t)}
=\frac{1}{\Gamma(1-\frac{\alpha}{2})}.
$$
\end{lemma}

\begin{thm}\label{thm:d-b arithmetic}
Let $\alpha=d- \b\in (0,1)$, where $\b$ is the constant in \eqref{eqn:MD} and $G$ is a set given as \eqref{eqn:SS1}
with $G_{0}$ being an open interval when 
$d=1$ or a bounded $C^{1,1}$ open set when $d\geq 2$.
Assume that $\{\ln(1/r_{j})\}_{j=1}^{N}$ is arithmetic with span $\rho$. 
Define $A=\sup_{z\in \R}s(z)$ and $B=\inf_{z\in\R}s(z)$, where $s(z)$ is from Theorem \ref{thm:BM}. 
\begin{enumerate}[(1)]
\item Let $g(t):=\int_{0}^{t^{-2/\alpha}}s(-\ln(t^{2/\alpha}v))v^{\alpha/2}\P(S_{1}^{(\alpha/2)}\in dv)$. Then, we have
\beq\label{eqn:critical asymp}
|G|-\tilde{Q}_{G}^{(\alpha)}(t)=tg(t) +o(t\ln(1/t)).
\eeq

\item\label{eqn:thm second part}  We have
\beq\label{eqn:critical oscillation}
\limsup_{t\to 0}\frac{g(t)}{\ln(1/t)} = \frac{A}{\Gamma(1-\frac{\alpha}{2})} \quad \text{ and } \quad
\liminf_{t\to 0}\frac{g(t)}{\ln(1/t)} = \frac{B}{\Gamma(1-\frac{\alpha}{2})}.
\eeq

\end{enumerate}
\end{thm}
\pf
Note that by the scaling property of $S_{t}^{(\alpha/2)}$ we have
\[
\begin{split}
&|G|-\tilde{Q}_{G}^{(\alpha)}(t)=\int_{0}^{\infty}\left(|G|-Q_{G}^{(2)}(u)\right)\P(S_{t}^{(\alpha/2)}\in du)
=\int_{0}^{\infty}\left(|G|-Q_{G}^{(2)}(t^{2/\alpha}v)\right)\P(S_{1}^{(\alpha/2)}\in dv)\nn\\
&= \int_{0}^{ t^{-2/\alpha}}\left(|G|-Q_{G}^{(2)}(t^{2/\alpha}v)\right)\P(S_{1}^{(\alpha/2)}\in dv) + 
\int_{ t^{-2/\alpha}}^{\infty}\left(|G|-Q_{G}^{(2)}(t^{2/\alpha}v)\right)\P(S_{1}^{(\alpha/2)}\in dv)\nn\\
&=\int_{0}^{t^{-2/\alpha}}\frac{|G|-Q_{G}^{(2)}(t^{2/\alpha} v)}{(t^{2/\alpha} v)^{\frac{d-\b}{2}}}
(t^{2/\alpha} v)^{\frac{d-\b}{2}}\P(S_{1}^{(\alpha/2)}\in dv) 
+\int_{t^{-2/\alpha}}^{\infty}\left(|G|-Q_{G}^{(2)}(t^{2/\alpha}v)\right)\P(S_{1}^{(\alpha/2)}\in dv).
\end{split}
\]
Hence, we have
\[
\begin{split}
|G|-\tilde{Q}_{G}^{(\alpha)}(t)-tg(t)
&=t\int_{0}^{t^{-2/\alpha}}\bigg(\frac{|G|-Q_{G}^{(2)}(t^{2/\alpha} v)}{(t^{2/\alpha} v)^{\frac{d-\b}{2}}}-s(-\ln(t^{2/\alpha}v))\bigg) 
v^{\frac{d-\b}{2}}\P(S_{1}^{(\alpha/2)}\in dv) \\
&\qquad \qquad 
+\int_{t^{-2/\alpha}}^{\infty}\left(|G|-Q_{G}^{(2)}(t^{2/\alpha}v)\right)\P(S_{1}^{(\alpha/2)}\in du).
\end{split}
\]

It follows from \cite[Equation (2.8)]{PS} we have 
\beq\label{eqn:limit small term}
\int_{t^{-2/\alpha}}^{\infty}\left(|G|-Q_{G}^{(2)}(t^{2/\alpha}v)\right)\P(S_{1}^{(\alpha/2)}\in du)
\leq c\int_{t^{-2/\alpha}}^{\infty}|G|u^{-1-\frac{\alpha}{2}}du
=o(t\ln(1/t)).
\eeq

In this case, by applying Theorem \ref{thm:BM} we have 
\beq\label{eqn:arithmetic expression}
\frac{|G|-Q_{G}^{(2)}(t^{2/\alpha} v)}{(t^{2/\alpha} v)^{\frac{d- \b}{2}}}=
s(-\ln(t^{2/\alpha}v)) +o(1)\quad \text{ as } t\to0.
\eeq
From \eqref{eqn:arithmetic expression} for any $\eps>0$ there exists $t_{0}(\eps)$ such that 
$$
\left|\frac{|G|-Q_{G}^{(2)}(t^{2/\alpha} v)}{(t^{2/\alpha} v)^{\frac{d-\b}{2}}}-s(-\ln(t^{2/\alpha}v))\right|<\eps
$$
for all $t\leq t_{0}$. Hence it follows from Lemma \ref{lemma:SS expectation} 
\[
\begin{split}
&\limsup_{t\to 0}\frac{t\int_{0}^{t^{-2/\alpha}}\Big(\frac{|G|-Q_{G}^{(2)}(t^{2/\alpha} v)}{(t^{2/\alpha} v)^{\frac{d-\b}{2}}}
-s(-\ln(t^{2/\alpha}v))\Big) v^{\frac{d-\b}{2}}\P(S_{1}^{(\alpha/2)}\in dv)}{t\ln(1/t)}\\
&\leq \eps\limsup_{t\to 0}\frac{\int_{0}^{t^{-2/\alpha}} v^{\alpha/2} \P(S_{1}^{(\alpha/2)}\in dv)}{\ln(1/t)}
\leq \frac{\eps}{\Gamma(1-\frac{\alpha}{2})}.
\end{split}
\]
This establishes \eqref{eqn:critical asymp}.

For \eqref{eqn:critical oscillation}, note that
\begin{equation}\label{eqn:d-b expression}
\begin{split}
|G|-\tilde{Q}_{G}^{(\alpha)}(t)\nn &
= \int_{0}^{ t^{-2/\alpha}}\frac{|G|-Q_{G}^{(2)}(t^{2/\alpha} v)}{(t^{2/\alpha} v)^{\frac{d-\b}{2}}}(t^{2/\alpha} v)^{\frac{d-\b}{2}}
\P(S_{1}^{(\alpha/2)}\in dv)\\
&\qquad \qquad 
 +\int_{ t^{-2/\alpha}}^{\infty}\left(|G|-Q_{G}^{(2)}(t^{2/\alpha}v)\right)\P(S_{1}^{(\alpha/2)}\in du).
\end{split}
\end{equation}
As in \eqref{eqn:limit small term} the second expression above is $o(t\ln(1/t))$ as $t\to 0$. 
For any $\eps>0$ it follows from \eqref{eqn:arithmetic expression} we have
 $\frac{|G|-Q_{G}^{(2)}(t^{2/\alpha} v)}{(t^{2/\alpha} v)^{\frac{d-\b}{2}}} <A+\eps$ 
for all sufficiently small $t$. 
This fact, together with \eqref{eqn:critical asymp} and Lemma \ref{lemma:SS expectation} gives 
$$
\limsup_{t\to0}\frac{g(t)}{\ln(1/t)}=\limsup_{t\to 0}\frac{|G|-\tilde{Q}_{G}^{(\alpha)}(t)}{t\ln(1/t)}\leq \frac{A+\eps}
{\Gamma(1-\frac{d-\b}{2})}.
$$
Since $\eps>0$ is arbitrary, we conclude that 
\beq\label{eqn:d-b limsup}
\limsup_{t\to0}\frac{g(t)}{\ln(1/t)}=\limsup_{t\to 0}\frac{|G|-\tilde{Q}_{G}^{(\alpha)}(t)}{t\ln(1/t)}\leq \frac{A}{\Gamma(1-\frac{d-\b}{2})}.
\eeq
For the lower bound, it follows from Theorem \ref{thm:BM} that for any $\eps>0$ there exists a sequence $t_{n}\to 0$ such that 
$$
\frac{|G|-Q_{G}^{(2)}(t_n^{2/\alpha} v)}{(t_n^{2/\alpha} v)^{\frac{d-\b}{2}}} \geq A-\eps.
$$
Hence, from Lemma \ref{lemma:SS expectation}, \eqref{eqn:limit small term},and \eqref{eqn:d-b expression} we have 
\begin{equation*} 
\begin{split}
\limsup_{n\to\infty}\frac{g(t_n)}{\ln(1/t_n)} &=
\limsup_{n\to \infty}\frac{|G|-\tilde{Q}_{G}^{(\alpha)}(t_n)}{t_{n}\ln(1/t_n)}\\
&\geq (A-\eps)\limsup_{n\to\infty}\int_{0}^{t_{n}^{-2/\alpha}}v^{\frac{\alpha}{2}}\P(S_{1}^{(\alpha/2)}\in dv)
\geq\frac{A-\eps}{\Gamma(1-\frac{d-\b}{2})}.
\end{split}
\end{equation*}
Since $\eps>0$ is arbitrary, we have 
\beq\label{eqn:d-b liminf}
\limsup_{t\to0}\frac{g(t)}{\ln(1/t)}\geq \frac{A}{\Gamma(1-\frac{d-\b}{2})}.
\eeq
Hence, the $\limsup$ version of \eqref{eqn:critical oscillation} follows from \eqref{eqn:d-b limsup} and \eqref{eqn:d-b liminf}, and 
the $\liminf$ version can be proved in the same way.
\qed

Here is the result for the non-arithmetic case.
The proof is very similar to the proof of Theorem \ref{thm:d-b arithmetic}, hence it will be omitted.
\begin{thm}\label{thm:d-b non-arithmetic}
Let $\alpha=d- \b\in (0,1)$, where $\b$ is the constant in \eqref{eqn:MD} and $G$ is a set given as \eqref{eqn:SS1}
with $G_{0}$ being an open interval when $d=1$ or a bounded $C^{1,1}$ open set when $d\geq 2$.
Assume that $\{\ln(1/r_{j}\}_{j=1}^{N}$ is non-arithmetic. 
Then, we have
$$
\lim_{t\to0}\frac{|G|-\tilde{Q}_{G}^{(d-\b)}(t)}{t\ln(1/t)}=\frac{1}{2 \Gamma(1-\frac{d-\b}{2}) 
\sum_{j=1}^{N}r_{j}^{\b}\ln(1/r_{j})}\int_{0}^{\infty}
\left(|G_{0}|-Q_{G_{0}}^{(2)}(u)\right)u^{-1-\frac{d-\b}{2}}du.
$$
\end{thm}

\subsection{The case of $\alpha \in (0,d-\b)$}
Now we handle the case when $\alpha\in (0,d-\b)$. 
We need the following simple lemma for the continuity of the map $t\to \tilde{Q}_{D}(t)$ which is proved 
in \cite[Lemma 3.11]{KP}.
\begin{lemma}\label{lemma:continuous}
For any open set $D$ with $|D|<\infty$, the map $t\to \tilde{Q}^{(\alpha)}_{D}(t)$ is continuous.
\end{lemma}

The following proposition is proved in \cite[Proposition 3.12]{KP}.
\begin{prop}\label{prop:alpha small}
Let $f$ be a bounded continuous function on $(0,\infty)$ such that $\displaystyle\lim_{x\downarrow 0}\frac{f(x)}{x^{\gamma}}$ 
exists as a finite number for some constant $\gamma>\frac{\alpha}{2}$. Then, we have 
$$
\lim_{t\downarrow 0}\int_{0}^{\infty}f(u)\frac{\P(S_{t}^{(\alpha/2)}\in du)}{t}= \frac{\alpha}
{2\Gamma(1-\frac{\alpha}{2})} \int_{0}^{\infty}f(u)\,u^{-1-\frac{\alpha}{2}}du.
$$
\end{prop}

\begin{thm}\label{thm:small alpha}
Let $\alpha\in (0,d- \b)$, where $\b$ is the constant in \eqref{eqn:MD} and $G$ is a set given as \eqref{eqn:SS1}
with $G_{0}$ being an open interval when 
$d=1$ or a bounded $C^{1,1}$ open set when $d\geq 2$.
Then, we have
$$
\lim_{t\to 0}\frac{|G|-\tilde{Q}_{G}^{(\alpha)}(t)}{t}= \frac{\alpha} {2\Gamma(1-\frac{\alpha}{2})} \int_{0}^{\infty}\big(|G|-Q_{G}^{(2)}(u)\big)\,
u^{-1-\frac{\alpha}{2}}du.
$$
\end{thm}
\pf
Note that we have 
$$
|G|-\tilde{Q}_{G}^{(\alpha)}(t)=\int_{0}^{\infty} \big(|G|-Q_{G}^{(2)}(u) \big)\P(S_{t}^{(\alpha/2)}\in du).
$$
It follows from Theorem \ref{thm:BM} that there exists constants $c_1$ such that
$$
|G|-Q_{G}^{(2)}(u) \leq  c_1u^{\frac{d-\b}{2}} \text{ for } u\leq1.
$$
Since $\alpha\in (0,d-\b)$, we can take $\gamma\in (\frac{\alpha}{2},\frac{d-\b}{2})$ and this implies 
$$
\lim_{u\to 0}\frac{|G|-Q_{G}^{(2)}(u)}{u^{\gamma}}=0.
$$
Now the conclusion of the theorem follows immediately from Proposition \ref{prop:alpha small}.
\qed

\medskip
\bigskip{\bf Acknowledgements}\,
The research of Y. Xiao is partially supported by NSF grant DMS-1855185.

\begin{singlespace}

\end{singlespace}
\end{doublespace}

\vskip 0.3truein

\noindent {\bf Hyunchul Park}

\noindent Department of Mathematics, State University of New York at New Paltz, NY 12561, USA

\noindent E-mail: \texttt{parkh@newpaltz.edu}

\bigskip

\noindent{\bf Yimin Xiao}

\noindent Department of Statistics and Probability, Michigan State University, East Lansing, MI 48824, USA

\noindent E-mail: \texttt{xiaoyimi@stt.msu.edu}


\begin{thebibliography}{99}

\bibitem{Val2016} 
L. Acu\~{n}a Valverde.
On the one dimensional spectral heat content for stable processes.
{\it J. Math. Anal. Appl.},  {\bf 441}  (2016),  11--24. 

\bibitem{V1} L. Acu\~{n}a Valverde.
Heat content for stable processes in domains of $\R^{d}$.
\textit{J. Geom. Anal.}, {\bf 27} (2017),  492--524.

\bibitem{B1994} M. van den Berg.
Heat content and Brownian motion for some regions with a fractal boundary.  
{\em Probab. Th. Rel. Fields}, {\bf 100} (1994), no. 4, 439--456. 

\bibitem{B2000} M. van den Berg.
Heat equation on the arithmetic von Koch snowflake. 
{\em Probab. Th. Rel. Fields}, {\bf 118} (2000), no. 1, 17--36. 


\bibitem{BD1989} 
M. van den Berg and E. B. Davies.
Heat flow out of regions in $\R^{m}$.
{\it Math. Z.}, {\bf 202} (1989), 463--482.


\bibitem{BG1998} M. van den Berg and P. B. Gilkey.
A comparison estimate for the heat equation with an application to the heat content of the 
$s$-adic von Koch snowflake. {\em Bull. London Math. Soc.},  {\bf 30} (1998), no. 4, 404--412.


\bibitem{FLV1995} 
J. Fleckinger, M. Levitin, and D. Vassiliev.
Heat equation on the triadic von Koch snowflake: asymptotic and numerical analysis. 
{\em Proc. London Math. Soc. (3)}, {\bf 71} (1995), no. 2, 372--396.

\bibitem{GPS19}
T. Grzywny, H. Park, and R. Song. Spectral heat content for L\'evy processes. 
{\em Math. Nachr.}  \textbf{292} (2019), 805--825.


\bibitem{KP}
K. Kobayashi and H. Park. 
Spectral heat contents for time-changed killed Brownian motions.
Preprint. \href{https://arxiv.org/abs/2007.05776}{arxiv.org/abs/2007.05776}.


\bibitem{Lapidus13}
M. L. Lapidus and M.  van Frankenhuijsen, {\it Fractal Geometry, Complex Dimensions 
and Zeta Functions. Geometry and spectra of fractal strings.} Second edition.  
Springer, New York, 2013.

\bibitem{Lapidus17}
M. L. Lapidus, G. Radunovi\'c,  and D. \u Zubrini\'c, {\it Fractal Zeta Functions and Fractal Drums. 
Higher-dimensional Theory of Complex Dimensions.}  Springer, Cham, 2017.


\bibitem{LV1996}
M. Levitean and D. Vassiliev.
Spectral asymptotics, renewal theorem, and the Berry conjecture for a class of fractals.
{\it Proc. London Math. Soc.} {\bf 72 (3)} (1996), no. 1, 188--214. 


\bibitem{PS}
H. Park and R. Song.
Small time asymptotics of spectral heat contents for subordinate killed Brownian 
motions related to isotropic $\alpha$-stable processes.
{\it Bul. London Math. Soc.},  {\bf 51}, Issue 2 (2019), 371--384.

\bibitem{PS21}
H. Park and R. Song.
Spectral heat content for $\alpha$-stable processes in $C^{1,1}$ open sets.
Preprint. \href{https://arxiv.org/abs/2007.02815}{arxiv.org/abs/2007.02815}. 

\end{thebibliography}
\end{document}